\documentclass{amsart}
\usepackage{amsfonts,amssymb,amsmath,amscd,amsthm}
\usepackage[utf8]{inputenc}
\usepackage{graphicx}
\usepackage{xcolor}
\usepackage{hyperref}
\usepackage{enumerate}
\usepackage{mathtools}

\title{Computability in Dynamical Systems}

\renewcommand{\epsilon}{\varepsilon}

\newcommand{\Per}{{\rm Per}}

\newcommand{\Or}{\mathcal{O}}
\newcommand{\C}{\mathbb{C}}

\newcommand{\cM}{\mathcal{M}}

\newcommand{\bR}{{\mathbb R}}
\newcommand{\bC}{{\mathbb C}}
\newcommand{\bZ}{{\mathbb Z}}
\newcommand{\bN}{{\mathbb N}}
\newcommand{\bQ}{{\mathbb Q}}

\newcommand{\cF}{{\mathcal F}}

\newcommand{\cG}{{\mathcal G}}
\newcommand{\cA}{{\mathcal A}}

\newcommand{\cL}{{\mathcal L}}

\newcommand{\RS}{\overline{{\mathbb C}}}

\newcommand{\eqdef}{{\vcentcolon=}}

\renewcommand{\epsilon}{\varepsilon}

\theoremstyle{definition}
\newtheorem{definition}{Definition}[section]
\newtheorem{theorem}[definition]{Theorem}
\newtheorem{proposition}[definition]{Proposition}
\newtheorem{corollary}[definition]{Corollary}

\newtheorem{example}[definition]{Example}
\begin{document}


\author{Michael Burr}
\author{Christian Wolf}

\address{Michael Burr, School of Mathematical and Statistical Sciences, Clemson University}
\email{burr2@clemson.edu}
\address{Christian Wolf, Department of Mathematics, The CUNY Graduate Center, NY, NY}
\email{cwolf@gc.cuny.edu}

\thanks{Burr was partially supported by grants from the National Science Foundation (CCF-1527193 and DMS-1913119).}

\thanks{Wolf was partially supported by  grants from  the Simons Foundation (\#637594 to Christian Wolf) and PSC-CUNY (TRADB-51-63715 to Christian Wolf).}

\maketitle

\begin{abstract}
In this paper we present an introduction to the area of computability in dynamical systems. This is a fairly new field which has received quite some attention in recent years.
One of the central questions in this area is if relevant dynamical objects can be algorithmically presented by a Turing machine.  
After providing an overview of the relevant objects we discuss recent results concerning the computability of the entropy for symbolic systems and the computability of Julia sets as well as their Brolin-Lyubich measures. 
\end{abstract}

\section{Introduction}

In \cite{Milnor:Remark}, Milnor poses the question of whether an algorithm exists to approximate the entropy of a dynamical system up to a specified error.  Such an algorithm does not exist in the most general setting, but this nonexistence leads to the following question:
\begin{center}
{\em For which classes of dynamical systems and invariants can these invariants be algorithmically approximated to any specified precision?}
\end{center}
In dynamics and the sciences, computational experiments are often used to gain insight into the behavior of dynamical systems.  For invariants which cannot be algorithmically approximated, the results of such computational experiments may be meaningless.  We approach the generalization of Milnor's question from the point of view of computability theory.

We begin with the computational model for computability theory.  From a theoretical standpoint, we may study a dynamical system $f:X\rightarrow X$ defined by an infinite amount of data.  From an algorithmic standpoint, however, only a finite amount of data may be provided as input.  Therefore, by necessity, we replace the dynamical system $f:X\rightarrow X$ by a convergent sequence of approximations $f_i:X_i\rightarrow X_i$, each of which may be described by a finite amount of data.  Milnor's question then leads to the question 
\begin{center}
{\em Do computations with a sequence of approximations $f_i:X_i\rightarrow X_i$ provide insight into the properties of the original dynamical system $f:X\rightarrow X$?}
\end{center}
We apply computability theory to three case studies from dynamical systems which yield both positive and negative results.


\subsection{Dynamical Systems}

A dynamical system consists of a state space $X$ (the set of attainable states) and a time evolution law
\[
(\Phi_t)_{t\in T}:X\to X,
\]
where $T$ denotes the set of times when the state of the system can be observed.  In other words, if the system is in state $x_0$ at time $t=0$, then it is in state $\Phi_t(x_0)$ at time $t$.  We focus on discrete dynamical systems where $T=\bN_0$.  We assume that $\Phi_0={\rm id}_X$ and that  $\Phi_{n+m}=\Phi_n\circ \Phi_m$ for all $n,m\in \bN$. Hence
\[
\Phi_n(x_0)=\underbrace{\Phi_1\circ \ldots \circ \Phi_1}_{n-{\rm times}}(x_0)=\vcentcolon\Phi_1^n(x_0).
\]
Consequently, such a dynamical system is given by the iteration of a map $f\vcentcolon=\Phi_1:X\rightarrow X$.

We consider discrete dynamical systems given by continuous maps $f:X\to X$ where $(X,d)$ is a compact metric space. 
Given $x\in X$ we call $\Or(x)\vcentcolon=\Or^+(x)\vcentcolon=\{x,f(x), f^2(x), \ldots\}$ the {\em (forward) orbit} of $x$. If $f$ is invertible, we analogously define the {\em backward orbit} $\Or^-(x)$ by iterating $f^{-1}$.  When $f^n(x)=x$ for some finite positive $n$, we say that $x$ is {\em periodic} and define the {\em period} of $x$ to be the cardinality of the orbit, i.e., a periodic point of period $n$ has orbit $\Or(x)=\{x,f(x),\dots,f^{n-1}(x)\}$, where each $f^i(x)$ is distinct for $0\leq i<n$ and $x=f^n(x)$.  We denote the set of {\em periodic points of $f$ with period $n$} by $\Per_n(f)$ and the set of {\em periodic points} of $f$ by $\Per(f)=\bigcup_n \Per_n(f)$.

In dynamics, it is a central goal to understand the asymptotic behavior of the orbits of $f:X\to X$ for a large class of points in $X$.  For instance, in ergodic theory, large often means a set of full measure for an $f$-invariant probability measure.

To motivate some of the challenges in studying families of dynamical systems computationally, the following example illustrates one of the pitfalls in the computational study of a single dynamical system.

\begin{example}
Consider the {\em $2x$-mod-$1$} dynamical system where $X=[0,1)$ and $f:X\rightarrow X$ takes $x$ to the fractional part of $2x$.  It is straightforward to show that some iterate of $x$ is periodic if and only if $x$ is rational.  Moreover, every dyadic rational is eventually mapped onto the fixed point $0$. Conversely, the closure of the orbit of an irrational point is either a Cantor set or  $X$ itself.
\end{example}

One challenge in studying this dynamical system computationally is that the orbits of rational points are not representative of a typical orbit, even though the rationals are dense in $X$.  Also, the most common number types on a computer are dyadic rationals, which exhibit even more atypical behavior.  This example illustrates that care must be taken with both the selection of representative points and how conclusions are drawn from computational experiments.

\subsection{Computability theory}
The main technique in computability theory is to replace a single mathematical object by a convergent sequence, and to use the convergent sequence to study the original object, see \cite{Brattka2008,RY1}.
\begin{definition}
An {\em oracle approximation} of a real number $\alpha$ is a function $\psi$ such that for any $n\in\mathbb{N}$, $\psi(n)$ is a rational number such that $|\alpha-\psi(n)|<2^{-n}$.  A real number $\alpha$ is {\em computable} if there exists a Turing machine which is an oracle approximation for $\alpha.$
\end{definition}
For most purposes, it is sufficient to think of a Turing machine as an algorithm.  Since every Turing machine is represented by a finite string of characters from a finite alphabet, we conclude that the computable numbers form a countable dense subset of the real numbers.  This conclusion immediately implies that noncomputable real numbers exist.

The error estimates in the definition of a computable real number are quite strong.  In some cases, there is a Turing machine that approximates $\alpha$ without an estimate on the convergence rate.
\begin{definition}
A real number $\alpha$ is {\em computable from above} (respectively, {\em from below}) if there exists a Turing machine $\psi$ such that on input $n$, $\psi(n)$ is a rational number such that the sequence $(\psi(n))_{n\in\mathbb{N}}$ is a decreasing (an increasing) sequence converging to $\alpha$.
\end{definition}
We observe that computability is equivalent to simultaneous computability from above and below.  The ideas from computability theory can be extended to functions as well as to many other settings.
\begin{definition}
A function $f:\bR\rightarrow\bR$ is a {\em computable function} if there exists a Turing machine $\phi$ with the following property: For any $\alpha\in\bR$ and oracle approximation $\psi$ for $\alpha$, on input $(\psi,n)$, $\phi(\psi,n)$ is a rational number such that $|f(\alpha)-\phi(\psi,n)|<2^{-n}$.
\end{definition}
The Turing machine $\phi$ must be able to accept every oracle approximation for any $\alpha$ (not merely Turing machines) and produce meaningful results.  Computable functions are necessarily continuous, and there are analogous notions of computability from above and below for computable functions.

We introduce the concept of a computable metric space and  refer the reader to \cite{HR} for more details.
\begin{definition}
Let $(X,d_X)$ be a separable metric space and let $\mathcal{S}_X=(s_i)_{i\in\mathbb{N}}$ be a dense sequence of points in $X$.  We say that $(X,d_X,\mathcal{S}_X)$ is a {\em computable metric space} if there is a Turing machine $\chi:\mathbb{N}^2\times\mathbb{N}\rightarrow\mathbb{Q}$ such that $|\chi(i,j,n)-d_X(s_i,s_j)|<2^{-n}$.
\end{definition}
The definition for a computable metric space is equivalent to the computability of the function $(i,j)\mapsto d_X(s_i,s_j)$.  Importantly the approximation is uniform, that is, the same Turing machine $\chi$ is used for every $i$ and $j$ in $\mathbb{N}$.  The definitions of computable numbers and functions generalize to computable metric spaces by replacing the output of Turing machines by an index $i\in\mathbb{N}$ of a point in $\mathcal{S}_X$.

\section{Computability of Entropy}
Arguably, the most important invariant of a dynamical system $f$ is its topological entropy $H_{\rm top}(f)$.  The topological entropy measures the complexity of the dynamical system $f$, and it is invariant under topological conjugation. Roughly, $H_{\rm top}(f)$ measures the exponential growth rate of the number of orbit segments that can be distinguished up to an error of $\epsilon$.  Specifically, for $n\in \bN$ and $\epsilon>0$, a finite set $F\subset X$ is {\em $(n,\epsilon)$-separated} if $\max_{k=0,\ldots,n-1} d(f^k(x),f^k(y))\geq \epsilon$ for distinct points $x,y\in F$. Let $F_n(\epsilon)$ denote the largest cardinality of an $(n,\epsilon)$-separated set. The {\em topological entropy} of $f$ is
\begin{equation}\label{defentropy}
  H_{\rm top}(f)\vcentcolon=  \lim_{\epsilon\to 0}\limsup_{n\to \infty}\frac{1}{n} \log F_n(\epsilon).
\end{equation}

From a computability standpoint, the limits in the formula for the topological entropy make its calculation difficult.  Even if $F_n(\epsilon)$ can be approximated for $n$ and $\varepsilon$, a computation cannot decide without additional theory whether $\epsilon$ is small enough and $n$ is large enough for $F_n(\epsilon)$ to approximate the limit.

We say that $f$ is {\em expansive} with {\em expansivity constant} $2\epsilon_0$ if for distinct $x,y\in X$, there is some $n$ such that $d(f^n(x),f^n(y))>2\epsilon_0$.  In this case, we can omit the limit on $\epsilon$ and compute the $\limsup$ for $\epsilon_0$.
Then sequence $(h_n)_n$ defined by 
\begin{equation}
    h_n\vcentcolon= \min_{k\leq n} \frac{1}{k}\log F_k(\epsilon_0) 
\end{equation}
is non-increasing and
converges 
to $H_{\rm top}(f)$.
\begin{corollary}\label{corupentropy}
Let $f:X\to X$ be an expansive dynamical system with expansivity constant $2\epsilon_0$ that is given by an oracle $\chi$. If there exists a Turing machine that takes $\chi$ and $n\in\bN$
as input and computes $F_n(\epsilon_0)$, then $H_{\rm top}(f)$ is upper semi-computable.
\end{corollary}

\subsection{Symbolic dynamical systems}
 
The class of symbolic systems form a rich family of discrete dynamical systems.  These systems provide good test cases for many theories as they can be studied in a hands-on way, but have enough variability to provide challenging examples.
We review the relevant material from symbolic dynamics, see \cite{MarcusLind1995}.

Let $\mathcal{A}=\{0,\dots,d-1\}$ be a finite set called the {\em alphabet}.  
Let $\cA^\ast$ denote the set of all finite words in $\cA$.
The {\em full shift} $\Sigma$ with alphabet $\mathcal{A}$ is the set of bi-infinite sequences $(x_k)_{k\in \bZ}$ of elements of $\cA$. 
With the {\em Tychonov product topology}, which coincides with the topology generated by the metric 
$$
d(x,y)\vcentcolon=2^{-\min\{|k| \;:\;  x_k\neq y_k\}},
$$
$\Sigma$ is a compact metrizable topological space. 

The  {\em (left) shift map} $f:\Sigma \to \Sigma $ is defined by $f(x)_k=x_{k+1}$, and we note that $f$ is a homeomorphism. We call a closed shift-invariant set $X\subset \Sigma$ a {\em shift space}, and  we say that $f:X\to X$ is a {\em subshift}.


For a fixed subshift $X$, let $\mathcal{L}(X,n)$ be the set of length-$n$ words which occur in some  element of  $X$.  Moreover, we call $\mathcal{L}=\bigcup_{n\in \bN} \mathcal{L}(X,n)$ the {\em language} of $X$. Languages and shift spaces are equivalent objects.

For a shift map $f:X\to X$, the entropy formula in Definition \eqref{defentropy} reduces to
\begin{equation}\label{eq:entropysymbolic}
H_{\rm top}(f|_X)\vcentcolon=\lim_{n\rightarrow\infty}\frac{1}{n}\log \#\mathcal{L}(X,n).
\end{equation}
Moreover, this limit converges from above.

An oracle $\chi$ of the language $\mathcal{L}(X)$ is a function such that $\chi(n)$ is a list of the words of the language of length $n$.
In light of the simplicity of the entropy of a subshift $X$ and Corollary \ref{corupentropy}, one might suspect that if the language of shift space is given by an oracle, then the entropy $H_{\rm top}(f|_X)$ is computable. This statement, however, is not true.

The following theorem uses the notion of {\em computability at a point}, which is a computable version of continuity at a point, see \cite{BurrWolf2018} for more details.

\begin{theorem}[{\cite{BDWY}}]\label{condcomputH(X)}
The generalized entropy map $X\mapsto H(f|_X)$ is computable at a nonempty shift space $X_0$ if and only if $X_0$ has zero entropy.  Moreover, there is one Turing machine that uniformly computes the topological entropy at all shift spaces with zero entropy.
\end{theorem}

This fact is proved by noting that the entropy is nonnegative, and, for any value of $\cL(X,n)$, there is a zero-entropy shift space $X'$ which has the same language as $X$ up to length $n$.  For example, $X'$ can be a union of orbits of points whose forward and backward tails are periodic.  Thus, given only the language, the only possible lower bound on the entropy is zero.  

\subsection{SFTs of finite type and Sofic shifts}
Subshifts of finite type and Sofic shifts are among the shifts for which the computability of their dynamics is best understood.  Both of these shifts can be described in terms of finite data.

Let $\mathcal{F}\subset\mathcal{A}^\ast$ be called a {\em forbidden set}. Then $X_\mathcal{F}$ is defined to be the largest shift space such that no point $x\in X_\mathcal{F}$ contains any of the words in $\mathcal{F}$. 

\begin{definition}
A shift space $X$ is a {\em subshift of finite type (SFT)} if  $X=X_\mathcal{F}$ for some finite set of words $\mathcal{F}$.
\end{definition}
The set $\mathcal{F}$ in this definition is not unique.
 A set $\mathcal{F}$ is a {\em minimal forbidden set} for an SFT $X$ if for every word in $\cF$, all of its proper subwords appear in the language of $X$.  It follows that there is only one minimal forbidden set for an SFT. 
 One less than the length of the largest word in the minimal forbidden set is called the {\em step size} of the SFT. Every SFT $X=X_{\mathcal{F}}$ is topologically conjugate to a one-step SFT $X'=X'_{\mathcal{F'}}$ on a possibly larger alphabet $\mathcal{A'}$. Moreover, given $\mathcal{F}$, the forbidden set $\mathcal{F'}$ can be computed by a Turing machine.
Thus, it suffices to study one-step SFTs. One-step SFTs can also be described via a {\em transition matrix}, as follows:
Let $A=(a_{i,j})$ be an $\cA\times\cA$
zero-one matrix. Then
\begin{equation}
X_A\vcentcolon=\{(x_k)_k: a_{x_i,x_{i+1}} =1\}
\end{equation}
is a one-step SFT. Moreover, every one-step SFT can be written this way. Additionally, the entropy $H(f|_{X_A})$ of the SFT $X_A$ is the logarithm of the spectral radius of the matrix $A$. Since there are algorithms that compute the spectral radius of $A$ at any given precision, we obtain the following result.

\begin{proposition}\label{prop:computabilitySFT:entropy}
There is a Turing machine that computes the topological entropy of any SFT given by either a finite forbidden set or a transition matrix. 
\end{proposition}


\begin{definition} Let $\mathcal{T}=(G,E,L)$ be a finite labeled directed graph, where $G$ is a graph with directed edge set $E$ and {\em labeling function} $L: E\rightarrow \mathcal{A}$, which assigns a {\em label} $L(e)$ from the finite alphabet $\mathcal{A}$ to each edge $e\in E$. Let $\xi=\cdots e_{-1}e_0e_1\cdots$ be a bi-infinite path on $G$
.  The {\em label} of the path $\xi$ is
$$L(\xi)\vcentcolon=\cdots L(e_{-1})L(e_0)L(e_1)\cdots\in \mathcal{A}^{\mathbb{Z}}.$$
The set of all bi-infinite labels of paths is denoted by  $X_{\mathcal{T}}=\{L(\xi): \xi\in G\}$.
A subset $X$ of the full shift is a {\em Sofic shift}  if $X=X_{\mathcal{T}}$ for a finite labeled graph $\mathcal{T}$.
\end{definition}

A shift space $X$ is Sofic if and only if $X$ is a (finite-to-one) factor of some SFT $Y$. In particular, every SFT is Sofic, but the converse is not true.
Moreover, there is a Turing machine that takes, as input, a finite labeled directed graph $\mathcal{T}$ and outputs  a transition matrix $A$ such that the Sofic shift $X_\mathcal{T}$ is a factor of the SFT $X_A$. Since finite-to-one factors preserve entropy, we obtain the following result.
\begin{proposition}
There exists a Turing machine that computes the topological entropy of any Sofic shift that is given by a finite labeled directed graph. 
\end{proposition}

\subsection{Invariant measures}\label{sec:invariantMeasures}

There is a natural connection between topological and measure-theoretic entropies, and we are interested in the computability properties of the measure-theoretic entropy.   A  probability measure $\mu$ is said to be {\em $f$-invariant} if $\mu(A)=\mu(f^{-1}(A))$ for all Borel sets $A$. 
For example, if $x\in \Per_n(f)$ then $\mu_x=\frac{1}{n}\sum_{k=0}^{n-1}\delta_{f^k(x)}$ is  the unique invariant measure supported on $\Or(x)$, where $\delta_{f^k(x)}$ denotes the Dirac measure on $f^{k}(x)$. 
Let $\cM_f$ denote the space of all $f$-invariant Borel probability measures  
on $X$ endowed with the weak$^\ast$ topology. With this topology, $\cM_f$ is a compact convex metrizable topological space. The extreme points in this convex space are denoted by $\cM_{f,E}$.  These measures are precisely the ergodic measures, where an $f$-invariant measure is {\em ergodic} if every $f$-invariant set has either measure zero or one.
\begin{example}
Let $X=[0,1)$ and $\alpha\in [0,\infty).$ Let $f:X\rightarrow X$ be a translation by $\alpha$, i.e., $f(x)$ is the fractional part of $x+\alpha$.  If $\alpha\in \bQ$,  then every point in $X$ is periodic and $\cM_{f,E}=\{\mu_x: x\in X\}\cup\{\lambda\}$, where $\lambda$ denotes the Lebesgue measure. On the other hand, if $\alpha\in \bR^+\setminus \bQ$ then
$\cM_f=\cM_{f,E}=\{\lambda\}.$
\end{example}
Let $H_\mu(f)$ denote the measure-theoretic entropy of $f$ with respect to $\mu$.  Roughly, $h_\mu(f)$ is similar to the topological entropy, but ignores sets of measure zero, see Section 4 for details.  The relationship between the topological and measure-theoretic entropies is the variational principle for the entropy, see Equation (\ref{varent}).


\subsection{General and coded shifts}
The computability properties of more general classes of shifts are more complicated.  One reason why these spaces are more challenging to study is that their descriptions may be given in terms of infinite amounts of data.  By Corollary \ref{corupentropy}, if the language is given by an oracle, then the topological entropy is upper semi-computable.  Thus, we are interested in identifying when the entropy is also lower semi-computable.

We approach this problem by using a sequence $\{X_m\}_{m\in\bN}$ of Sofic shifts which ``exhausts" $X$ in the sense that the entropy of these $X_m$'s approaches that of $X$ from below.  The challenge with this approach is that the generalized entropy function $X\mapsto H_{\rm top}(f|_X)$ might not be continuous.  Whenever such an approximation exists, one calls such a shift an {\em almost Sofic shift}, see \cite{Petersen1986}.  In most cases, the language itself is not enough to compute such an approximation, we, therefore, focus on the case where $X$ is generated from concatenations of a countable set of finite words, that is $X$ is a coded shift.

To emphasize the generality of the class of coded shifts, we note that every shift space $X$ can be expressed as the closure of the set of
bi-infinite paths on a countable directed labeled graph, called  a {\em representation} of $X$, see \cite{MarcusLind1995}.
Coded shifts are precisely those shift spaces which have a representation by an irreducible countable directed labeled graph.
Many classes of shifts spaces are coded shifts, including Sofic shifts, S-gap, and Beta-shifts.

\begin{definition}
A two-sided shift space $X$ over a finite alphabet $\cA=\{0,\dots,d-1\}$ is {\em coded} if there exists a {\em generating set} $\cG=\{g_i\}_{i\in\bN}\subset\cA^\ast$, such that $X=X(\cG)$ is the smallest shift space that contains all bi-infinite concatenations of generators, i.e.,
\[ X_{\rm seq}\eqdef\{ \cdots g_{i_{-2}}g_{i_{-1}}g_{i_0}g_{i_1}g_{i_2}\cdots: i_j\in \bN, j\in \bZ\}.\]
In other words, $X$ is the topological closure of $X_{\text{seq}}$, for more details, see \cite{BHcoded,MarcusLind1995}.
\end{definition}

Let $X_m\vcentcolon=X(\{g_1,\dots,g_m\})$ be the coded shift generated by the first $m$ generators in $\cG$.  Since each $X_m$ is a Sofic shift, its entropy is computable.  We observe that $\{X_m\}_{m\in\bN}$ is an increasing sequence of shift spaces in $X_{\rm seq}$.  In order to establish computability results for the entropy of coded shifts, we study when the entropy of $\{X_m\}_{m\in\bN}$ converges to the entropy of $X$.  For this convergence to hold, we show that the entropy on $X$ cannot be ``concentrated" on $X\setminus \bigcup_m X_m$.

To make these notions more precise, 
we introduce the following two definitions:
\begin{definition}\label{def:FSP}[{\cite{BDWY}}]
We say that $X$ has {\em full sequential entropy} if 
$$
H_{\rm top}(f|_X)=\sup_{\mu\in\cM}\{H_\mu(f):\mu(X_{\rm seq})=1\}.
$$
We say that $X$ has {\em strict full sequential entropy} if 
\begin{multline*}
\sup_{\mu\in \cM}\{H_\mu(f): \mu(X_{\rm seq})=1\}>
\sup_{\mu\in \cM}\{H_\mu(f): \mu(X\setminus X_{\rm seq})=1\}.
\end{multline*}
\end{definition}

We observe that strict full sequential entropy implies full sequential entropy.  

\begin{definition}[cf \cite{Pavlov:UniqueRepresentability}]
A {\em unique representation} of a coded shift $X$ is a generating set $\cG$ of $X$ such that each $x\in X_{\rm seq}(\cG)$ can be uniquely written as an infinite concatenation of elements of $\cG$.
\end{definition} 
The manuscript \cite{BPR21} includes a proof that all coded shifts have a unique representation.

\begin{theorem}[{\cite{BDWY}}]\label{thm:computabilityentropycoded}
Let $X$ be a coded shift with unique representation $\cG$.
\begin{enumerate}[\rm(i)]
\item If  $X$ has strict full sequential entropy then $H_{\rm top}(f|_X)=\lim_{m\to\infty} H_{\rm top}(f|_{X_m})$. 
\item 
Suppose  $X$ has strict full sequential entropy. 
There exists a Turing machine, independent of $X$, that computes $H_{\rm top}(f|_X)$ from oracles for the language and generating set $\cG$ of $X$. In particular, $H_{\rm top}(f|_X)$ is computable; and 
\item If $X$ does not have full sequential entropy, then the entropy is not continuous at $X$ and hence is not computable at $X$.
\end{enumerate}
\end{theorem}

Since the entropy is not continuous as a function of the shift space, the main challenge in the proof of Statement (ii) of this theorem is to show that the sequence of entropies  for $X_m$  converges to that of $X$.  In the proof of Statement (iii), it is shown that the entropy for subshifts $X_m$ does not converge to the entropy of $f$ on $X$. We refer to \cite{BDWY, Petersen1986} for examples of coded shifts without full sequential entropy. In particular, there are coded shifts for which the topological entropy is not computable.


\subsection{Classification of entropies}
An alternative to computing the entropy of a dynamical system is to classify the computability properties of entropies of dynamical systems.  For instance, by Proposition \ref{prop:computabilitySFT:entropy}, the entropies of all SFTs are not only computable, but also are rational logarithms of Perron numbers, where a {\em Perron number} is a real zero of an integer polynomial such that the zero is greater than one and is the zero of largest magnitude.  The converse of this statement also holds, i.e., every rational logarithm of a Perron number is the entropy of some mixing SFT, see~\cite{LindClassification}.

A natural generalization of one-dimensional SFTs are multidimensional SFTs which are finite type  $\bZ^d$-actions.  In this case, the set of possible entropies coincides with the set of positive upper computable numbers \cite{hochman2010characterization}. 
Since there are upper computable numbers which are not computable, there exist symbolic systems whose entropy is not a computable number.
On
the other hand, the entropy of strongly irreducible higher dimensional SFTs is computable.

In more generality, for computable dynamical systems, see the definition in \cite{GHRS}, the entropy is always a supremum of a computable sequence of upper-computable numbers \cite{GHRS}.  Moreover, for every positive such number, there is a dynamical system over $\{0,1\}^\bN$ which has this entropy and for which the time evolution function is computable.


\section{Computability of Julia sets}
Computer generated pictures of the Mandelbrot set and Julia sets can be quite striking, but computability theory raises the question of whether the underlying sample used to generate these sets accurately represents their true shape.  Julia sets are prototypes of fractals consisting of the points with unpredictable dynamics.  Since the dynamical behavior of nearby points can differ, computations with sample points may not correctly capture the dynamical behavior.

Let $\RS=\C\cup\{\infty\}$ be the Riemann Sphere endowed with the spherical metric $d_{\RS}$.  Let 
$$f_c(z)=z^d+c_{d-2}z^{d-2}+\dots+c_0$$ be a polynomial of degree $d$, where $c=(c_{d-2},\dots,c_0)\in \C^{d-1}$.  We denote the space of polynomials of this form by $\operatorname{Pol}(d)$.

Fatou and Julia introduced a partition of $\RS$ into two completely invariant sets, $F_f$ and $J_f$, called the {\em Fatou} and {\em Julia} sets, respectively.  The Fatou set $F_f$ is the set of points $z\in \RS$ where the family of iterates of $f$ is equicontinuous in a neighborhood of $z$.  The Julia set $J_f$ is the complement of $F_f$ in the Riemann sphere, and thus, the Fatou set is open, while the Julia set is compact.  In particular, $z\in F_f$ if and only if for all $\epsilon>0$
there exists $\delta>0$ such that $d_{\RS} (z,w)<\delta$ implies $d_{\RS}(f^n(z),f^n(w))<\epsilon$ for all $n\in \bN$. In other words, the Fatou set is the set of points whose orbits are Lyapunov stable. The Julia set is the set of points whose orbits have sensitive dependence on the initial conditions, i.e., {\em chaotic dynamics}.  

When $d=2$, the polynomials are of the form $f_c(z)=z^2+c$. The Julia set of a quadratic polynomial is either connected or totally disconnected. The connectivity-locus
\[
\cM=\{c\in \C: J_c\,\, \text{is connected}\}
\]
is the famous {\em Mandelbrot set}. 
The {\em filled-in Julia set} of $f$ is defined by
\[
K_f=\{z\in \bC: \{f^n(z): n\in \bN_0\}\,\, \text{is bounded}\}.
\]
The iterates of the points in $U_f=\RS\setminus K_f$ converge to $\infty$.  Moreover,
\[
J_f=\partial K_f=\partial U_f\,\, \text{and}\,\, F_f=U_f\cup \operatorname{int} K_f.
\]

The sensitivity to the initial conditions makes computing the Julia set particularly challenging.

\begin{figure*}
\begin{center}
\includegraphics[height=1.7in]{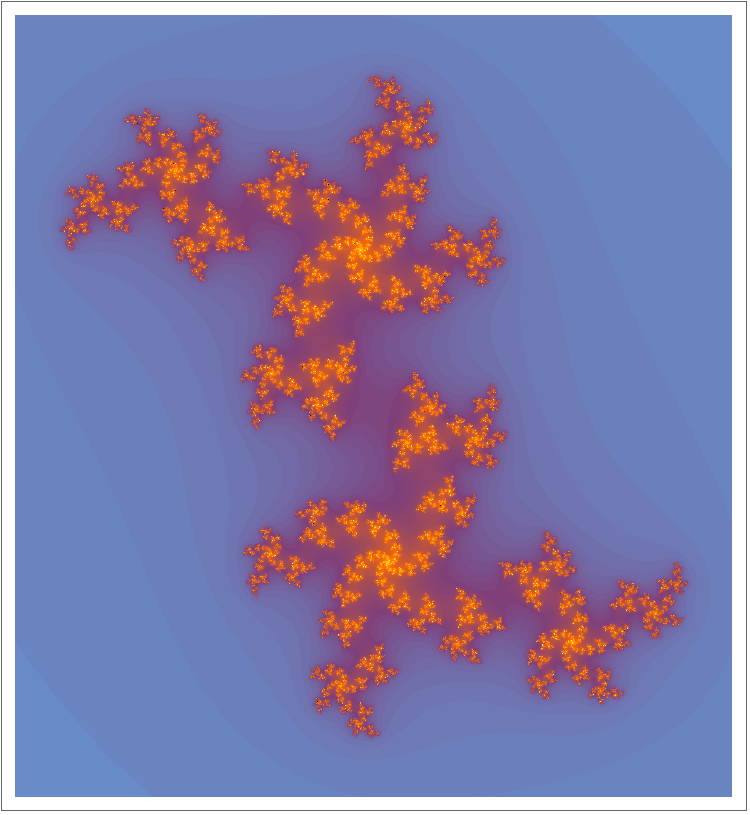}
\includegraphics[height=1.7in]{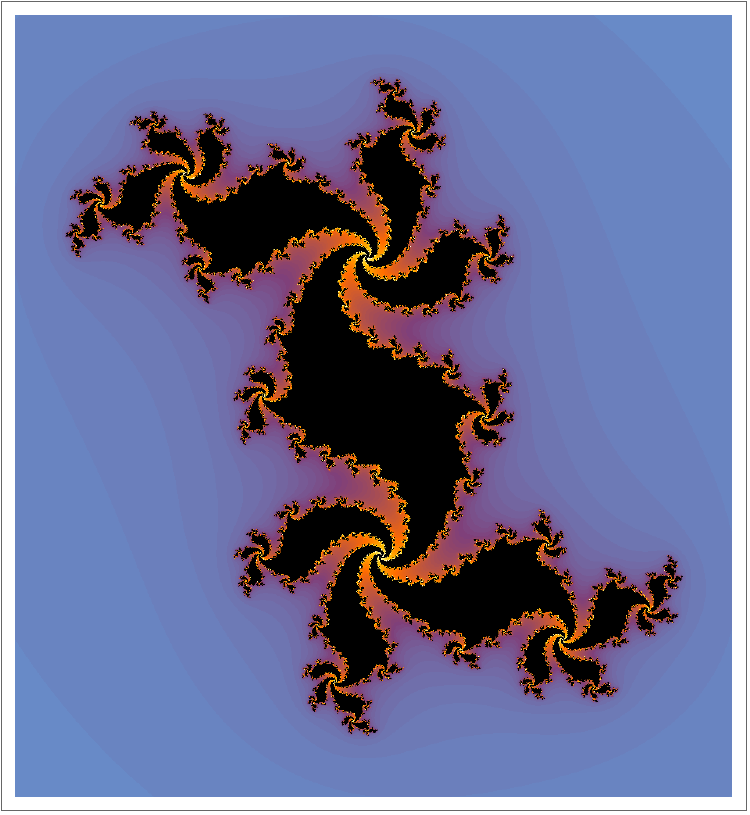}
\includegraphics[height=1.7in]{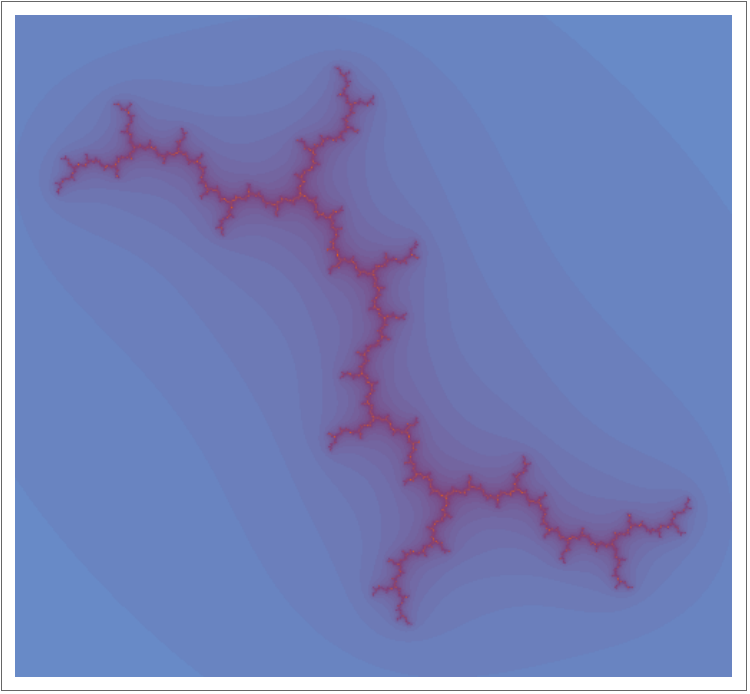}
\includegraphics[height=1.7in]{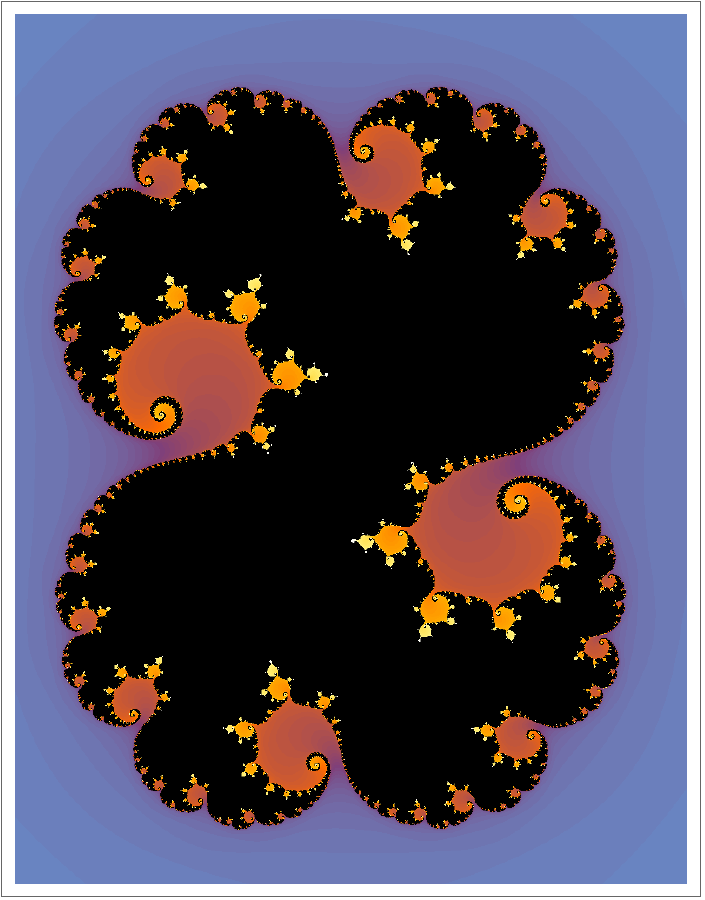}
\end{center}
\caption{Examples of (filled-in) Julia sets for $z^2+c$: (a) $c\not\in \cM$, totally disconnected Julia set (b) $c\in  {\rm int} \cM$, $f$ is hyperbolic   (c) $c\in \partial\cM$, Julia set is a dendrite (d) $c\in \partial \cM$, $0$ is parabolic fixed point. In all these examples the Julia set is computable. These images were created with Mathematica.}
\end{figure*}

\subsection{Computability of sets}

\begin{definition}
Let $S\subset\bR^m$ be an open set.  $S$ is a {\em recursively open set} (also called a semi-decidable set or a lower-computable set) if there exists a Turing machine $\psi$ that produces a (possibly infinite) sequence of pairs $(z_i,n_i)$ so that $z_i\in\bQ^m$ is a rational vector and $n_i\in\bZ$ so that $S=\cup_i B(z_i,2^{-n_i})$
\end{definition}

There is a Turing machine, whose input is $\psi$, that terminates if $s\in S$, but may fail to terminate if $s\not\in S$.

\begin{definition}
Let $S, T\subset \bR^m$  be compact sets.  The {\em Hausdorff distance} between $S$ and $T$ is 
$$d_H(S,T)\vcentcolon=\max\{\max_{s\in S}\min_{t\in T}d(s,t),\max_{t\in T}\min_{s\in S}d(s,t)\}.$$
\end{definition}

\begin{definition}
Let $S\subset \bR^m$ be compact.  $S$ is {\em computable} if there exists a Turing machine $\psi$ such that on input $n$, $\psi(n)$ is a finite collection of pairs ${(q_i,l_i)}$ where $q_i\in\bQ^m$ and $l_i\in\bZ$ such that $d_H(\cup_i \overline{B}(q_{k_i},2^{-l_i}),S)<2^{-n}$
\end{definition}

There is a Turing machine, whose input is $\psi$, that terminates if $s\not\in S$, but may fail to terminate if $s\in S$.

\subsection{Properties of periodic points}

 By a theorem of Fatou and Julia, the Julia set is the closure of the set of repelling periodic orbits. A periodic point $\alpha\in\Per_k(f)$ is {\em repelling} (respectively, {\em attracting}) if $|(f^k)'(\alpha)|>1$ (if $|(f^k)'(\alpha)|<1$).  A periodic point $\alpha\in\Per_k(f)$ is {\em neutral} if $|(f^k)'(\alpha)|=1$. 
 
For a neutral periodic point, we replace $f$ by its $k$-th iterate so that $\alpha$ is a fixed point.  We write $f'(\alpha)=e^{2\pi i \theta}.$ We say that $\alpha$ is {\em parabolic} if   $f'(\alpha)$ is a root of unity, i.e., $\theta\in \bQ$. In this case, $\alpha\in J$ and the dynamics near $\alpha$ is given by the Leau-Fatou Flower theorem, i.e., alternating attracting and repelling petals \cite{Milcomplexdynamicsbook}. 


We say that $\alpha$ is a {\em Siegel point} if the dynamics near $\alpha$ is conjugate to the linear map $L_\alpha(z)= f'(\alpha)z$. Then $\alpha\in F_f$  and the connected component $U\subset F_f$ containing $\alpha$ is a {\em Siegel disk}, i.e., there exists a conformal map $\varphi$ from $U$ to the unit disk $D$ such that $f(z)=\varphi^{-1}\circ L_\alpha \circ \varphi(z)$ for all $z\in U$.


By Bryuno's Theorem, $\alpha$ is Siegel if $\theta$ is a {\em Bryuno number}, that is, 
$$\sum_n \frac{\log (q_{n+1})}{q_n}<\infty,$$
where $q_n$ is the denominator of the $n$-th continued fraction expansion of $\theta$. Moreover, the Bryuno numbers have full Lebesgue measure.  We denote the set of degree $d$ polynomials with a Siegel disk by ${\rm Pol}_{\rm Sie}(d)$.

\subsection{Noncomputability of Julia sets}

The map $c\mapsto J_c$ is not continuous at $c$ if $f_c$ has either a parabolic or a Siegel point, see \cite{BY}. In particular, the map $c\mapsto J_c$ is not computable on ${\rm Pol}(d)$. Furthermore, for each non-computable real number $c$, the Julia set   $ J_c$ is not computable without oracle access to $c$. 

Braverman and Yampolsky show that $c\mapsto J_c$ is computable on ${\rm Pol}(d)\setminus {\rm Pol}_{\rm Sie}(d)$ by one Turing machine once some non-uniform information about each parabolic point $\alpha$ of $f_c$ is provided.
This shows that the computability of the Juila set is problematic only for polynomials with Siegel disks.  For example, consider the quadratic family
\[
f_\theta(z)=z^2+ e^{2\pi i\theta} z, \quad \theta \in [0,1),
\]
which has a neutral fixed point at $0$. If  $\theta$ is Bryuno, then $f_\theta$ has a Siegel disk centered at $0$. In particular, $f_\theta$ has a Siegel disk at the origin for $\theta$-values in a set of full Lebesgue measure. It is natural to conjecture that the non-computability of Siegel Julia sets $J_\theta$  stems from the non-computability of $\theta$ -- surprisingly, this is not the case. 

Braverman and Yampolsky also construct an explicit algorithm that computes $\theta$-values with a Siegel disk at the origin. In particular, there exist computable $\theta$ values with a noncomputable Julia set $J_\theta$. The idea of this proof is to show that the Julia set $J_\theta$ is computable if and only if the conformal radius $r(\theta)$ of its Siegel disk is computable.  In addition, they show that there are computable $\theta$'s such that $r(\theta)$ is not computable.  This shows that it is impossible to compute the Julia set from $\theta$ alone since then $r(\theta)$ could be computed from $\theta$, which is not possible.

\subsection{Computability of filled Julia sets}

For the filled-in Julia set the situation is better.
\begin{theorem}[\cite{BY}]\label{theKcomp}
For any $c$, there is a Turing machine that produces a $2^{-n}$-approximation of the filled-in Julia set $K_c$ from an oracle for $c$ and $n\in\bN$.
\end{theorem}
The proof of Theorem \ref{theKcomp} describes a pair of Turing machines, one which terminates on points further than $2^{-n}$ from the filled-in Julia set, and a second which terminates on points which are near the filled-in Julia set.  The second Turing machine classifies the behavior of points, e.g., near the Julia set or an iterate in the interior of a Siegel disk.

This result does not imply that $c\mapsto K_c$ is computable on ${\rm Pol}(d)$ since there would need to be {\it one} Turing machine that computes $K_c$ for {\it all} $c\in {\rm Pol}(d)$.

\medskip

There are other recent results that discuss the computability of dynamically defined sets, e.g., the computability of Lorenz attractors \cite{GRZ} and the computability of generalized rotation sets \cite{BSW}.










\section{Computability of measures of maximal entropy}
Let $f:X\to X$ be a continuous map on a compact metric space $(X,d)$.  As in Section \ref{sec:invariantMeasures}, there is the notion of the measure-theoretic entropy $H_\mu(f)$ of an invariant measure $\mu$, see Section \ref{sec:MaximalEntropy} for the definition.  A measure which maximizes $H_\mu(f)$ is called a {\em measure of maximal entropy.}  

Roughly, $\mu$ is a measure of maximal entropy if the set of its typical points carries the full topological complexity in terms of entropy. A dynamical system might not have a measure of maximal entropy, or, when such a measure exists, it might not be unique.

For several classes of dynamical systems, there exists a unique measure of maximal entropy, e.g.,  transitive SFTs and sofic shifts.  If $f:\RS\to\RS$ is a polynomial of degree $d\geq 2$, then $f$ has a unique measure of maximal entropy, called the {\em Brolin-Lyubich (B-L) measure} $\lambda_f$. The measure of maximal entropy $\lambda_f$ is fully supported on the Julia set $J_f$.

One might suspect that since there are non-computable Julia sets, the B-L measure is also not computable.  Surprisingly, this is not the case, and the B-L measure is computable.

\subsection{Measure-theoretic entropy}\label{sec:MaximalEntropy}
Using a result of Katok and Brin, we provide an intuitive description of $H_\mu(f)$ as a measure-theoretic notion of complexity.

For $n\in\bN$, define the Bowen $d_n$-metric on $X$ by
\[
d_n(x,y)=\max_{k=0,\dots,n-1} d(f^k(x),f^k(y)).
\]
We observe that $d_n$ induces the same topology on $X$ as $d$. Let $B_n(x,\epsilon)$ denote the Bowen-ball of radius $\epsilon$ with center $x$ in the $d_n$-metric. 

\begin{theorem}[Katok-Brin]
Let $\mu\in \cM_f$. Then 
\[
H_\mu(x)\vcentcolon=\lim_{\epsilon\to 0}\lim_{n\to\infty}- \frac{\log(\mu(B_n(x,\epsilon))}{n}
\]
exists for $\mu$ for almost every $x\in X$ and $H_\mu(\cdot)$ is $\mu$-measurable. Moreover, if $\mu$ is ergodic, then $H_\mu(\cdot)$ is constant almost everywhere. 
\end{theorem}

The quantity $H_\mu(x)$ is called the {\em point-wise entropy} of $\mu$ at $x$, and we define
$$
H_\mu(f)=\int_X H_\mu(x)\, d\mu(x).
$$
For $H_\mu(x)$ to be positive, the measure of the Bowen balls must decay exponentially to zero in $n$. Since Bowen balls are sets of points that stay close together under iteration, one infers that $H_\mu(f)>0$ implies that $\mu$ must be ``concentrated" on sets whose points are pushed apart under iteration. Thus, $\mu$ is concentrated on sets where a small change in the initial condition results in different dynamical behavior.

Bowen also defines the topology entropy $H_{\rm top}(Y)$ for arbitrary sets $Y\subset X$. If $\mu$ is ergodic, then
$H_\mu(f)=H_{\rm top}({\mathbb G}_\mu),$
where ${\mathbb G}_\mu$ denotes the set of $\mu$-generic points of $\mu$, where a point $x$ is {\em $\mu$-generic}
if its orbit distribution coincides with the measure $\mu$, i.e., 
\[
\lim_{n\to \infty}  \frac{1}{n}\sum_{k=0}^{n-1}\delta_{f^k(x)} = \mu,
\]
in the weak$^\ast$ topology.  As a consequence of Birkhoff's Ergodic theorem, for every ergodic measure $\mu$, the set of $\mu$-generic points has full $\mu$-measure. 
The topological and measure-theoretic entropies are related by the {\em variational principle for the entropy}:
\begin{equation}\label{varent}
    H_{\rm top}(f)=\sup_{\mu\in \cM_f} H_\mu(f).
\end{equation}

\subsection{Properties of B-L measures}

 The B-L measure $\lambda$ is characterized by the following equivalent properties, see e.g. \cite{BBRY}:
\begin{enumerate}
    \item $H_\lambda(f)=\log d=H_{\rm top}(f);$
    \item $\lambda(f(A))=d\cdot \lambda(A)$ for all measurable sets $A$ on which $f$ is injective;
    \item $\lim_{n\to\infty} \frac{1}{d^n} \sum_{w\in f^{-n}(z)} \delta_w =\lambda$ for all but finitely many $z\in\bC$ in the weak$^\ast$ topology;
    \item $\lambda$ is the harmonic measure  at infinity of the filled-in Julia set $K_f$.
\end{enumerate}

\subsection{Computability of measures}

Let $\cM(\RS)$ denote that space of Borel probability measures on $\RS$ endowed with the weak$^\ast$ topology. We endow $\cM(\RS)$
with the Wasserstein-Kantorovich
metric defined by
\[
W_1(\mu,\nu)\vcentcolon=\sup_{\phi\in 1-{\rm Lip}(\RS)} \left| \int \phi\, d\mu -\int \phi\, d\nu \right|,
\]
where $1{\rm -Lip}(\RS)$ is the space of real-valued Lipschitz functions on $\RS$ with Lipschitz constant $1$. Let ${\mathbb D}$ denote the set of probability measures on $\RS$ that are finite rational convex combinations of Dirac measures supported on points in $\bC$ with rational coordinates. Then $(\cM(\RS), W_1, {\mathbb D})$ is a computable metric space.

\subsection{Computability of B-L measures}

We address the computability of the map ${\rm Pol}(d)\ni c\mapsto \lambda_c$ from polynomials to B-L measures.

\begin{definition}
The function $c\mapsto \lambda_c$ is computable if there exists a Turing machine $\phi$ with the following property: For any $c\in {\rm Pol}(d)$ and any oracle approximation of $\psi$ of $c$, on input $\psi$ and $n$, $\phi(\psi,n)$ is a measure in ${\mathbb D}$ with $W_1(\phi(\psi,n),\lambda_c)<2^{-n}$.
\end{definition}

\begin{theorem}[{\cite{BBRY}}]\label{thmmmerat}
Let $d\geq 2$. Then $c\mapsto \lambda_c$ is computable.
\end{theorem}
The proof of Theorem \ref{thmmmerat} follows from general statements about the computability of measures in $\cM(\RS)$ and  from the characterization of the B-L measure in Property 2 of Section 4.2.

We note that Theorem \ref{thmmmerat} is surprising for parameters with  non-computable Julia sets. To see this consider a (non-exceptional) point $z_0\in\bC$. For any given $n\in \bN$, the pre-images $f_c^{-n}(z_0)$ can be calculated at any pre-described accuracy. Moreover, 
\begin{equation}\label{conv1}
    f_c^{-n}(z_0)\quad \longrightarrow \quad J_c
\end{equation}
in the Hausdorff metric. If $c$ is a parameter for which $J_c$ is not computable, then  the accuracy of the convergence in Equation \eqref{conv1} can not be controlled by a Turing machine, but the convergence of
\begin{equation}\label{conv2}
    \frac{1}{d^n}\sum_{w\in f_c^{-n}(z_0)} \delta_w \quad \longrightarrow \quad \lambda_c
\end{equation}
in the weak$^\ast$ topology can be computed at any prescribed accuracy.
Since the support of the B-L measure is the entire Julia set, we interpret the convergence in Equation~\eqref{conv2} as covering most of the Julia set at a prescribed rate which is sufficient to compute the   B-L measure.

\medskip

There are additional results in the literature concerning the computability and non-computability of other dynamically relevant measures including physical measures, e.g. \cite{GHR}.


\bibliographystyle{plain}
\bibliography{bibliography}

\end{document}